\documentclass[a4paper,12pt]{article}
\usepackage{amssymb,amsfonts,amsmath}
\usepackage[english]{babel}
\usepackage{latexsym}
\usepackage{epsfig}

\newcommand{\RR}{{\mathbb{R}}}
\newcommand{\ZZ}{{\mathbb{Z}}}

\newcommand{\QQ}{{\mathbb{Q}}}
\newcommand{\NN}{{\mathbb{N}}}

\title{Dissecting brick into bars.}
\author{Ivan Feshchenko, Danylo Radchenko, Lev Radzivilovsky, Maksym Tantsiura}

\begin{document}

\maketitle

\begin{abstract}
An $N$-dimensional parallelepiped will be called a bar if and only 
if there are no more than $k$ different numbers among the lengths of its sides
(the definition of bar depends on $k$). We prove that a parallelepiped can be 
dissected into finite number of bars iff the lengths of sides of the 
parallelepiped span a linear space of dimension no more than $k$ over $\QQ$.
This extends and generalizes a well-known theorem of Max Dehn about
partition of rectangles into squares. Several other results about 
dissections of parallelepipeds are obtained.
\end{abstract}


\section{Introduction} \label{sec:intro}
The following well-known result was proved by Dehn in 1903.

\smallskip
\textbf{Theorem 1 (Dehn).}  A rectangle can be dissected into squares iff the ratio 
of the width and the height is rational. 

\smallskip
Of course, one can talk about dissections of parallelograms into rhombs 
instead.

In this paper we investigate possible multidimensional generalizations of
theorem 1. 

Original Dehn's proof [7] was complicated.
Since 1903, several proofs were proposed for that theorem.
In 1940 Brooks, Smith, Stone and Tutte [5] found a surprising proof, 
in which they transform the question into a question about 
electric circuits.

In 1950's Hadwiger[10] found probably the shortest approach, using Hamel function [11]
to construct additive functions over rectangles. Similar proof was later 
published by Pokrovskii [19].
As Boltianskii [3],[4] pointed out, use of Axiom of Choice is unnecessary in those proofs. 

Another proof, based on deformations, which is not widely known,
belongs to A. Kanel-Belov (private communication). 

A nice survey, along with some other related theorems, can be found in 
the paper of Freiling and Rinne [8].

We shall reproduce Hadwiger/Pokrovskii proof because we shall use 
similar ideas in more general situations and proof of Kanel-Belov
because it is not widely known.

An obvious way to extend this theorem, which was noticed by several 
authors, is to ask when is it possible to cut 3-dimensional cuboid 
into cubes. Since it generates a partition of every face into squares,
it follows easily from Dehn's theorem that the ratio between each 
two sides of cuboid is rational.

An non-obvious extension of the theorem, which was a starting point
of our investigation, is dissecting parallelepipeds into bars.

\smallskip
\textbf{Definition.} A bar in $\RR^3$ is a parallelepiped which has no more 
than 2 different side lengths. 

Simply speaking, it is a box of type $a \times a \times b$.

\smallskip
\textbf{Theorem 2.} In 3-dimensional space, a parallelepiped with sides $a,b,c$ 
can be divided into bars iff there is a non-trivial linear combination 
of $a,b,c$ with integer coefficients which is $0$.
In other words, the condition is that a linear space over $\QQ$
spanned by $a,b,c$ is of dimension no more than 2.
\smallskip

One might try to generalize the notion of bar for 4-dimensional space in three ways: 
either cuboid of type $a \times a \times a \times b$ (three sides in different directions, 
doesn't matter which, are equal) or cuboid of type $a \times a \times b \times b$ or both.
Here we show 2 theorems, which give some intuition, why we should take both types,
to get a theorem similar to theorems 1 and 2:

\smallskip
\textbf{Theorem 3.} In 4-dimensional space, a parallelepiped of type $a \times a \times b \times b$ 
can be divided into parallelepipeds of type $x \times x \times x \times y$,  
iff there is a non-trivial linear combination of $a,b$ with integer coefficients 
which is $0$.

\smallskip
\textbf{Theorem 4.} In 4-dimensional space, a parallelepiped of type $a \times a \times a \times b$ 
can be divided into parallelepipeds of type $x \times x \times y \times y$ 
iff there is a non-trivial linear combination of $a,b$ with integer coefficients 
which is $0$.

\medskip
So, theorems 3 and 4 suggest that to extend theorem 2, we should use the following 
definition:

\smallskip
\textbf{Definition.} A $k$-bar in $n$-dimensional space is a cuboid with no more than $k$ different side lengths. 

\smallskip
\textbf{Theorem 5.} In $n$-dimensional space, a cuboid with sides $a_1,a_2,...,a_n$ can be divided into $k$-bars 
iff the dimension of $\QQ$-linear space, spanned by $a_1,a_2,...,a_n$ is not greater then $k$.
\smallskip

There is a well known fact, which looks somewhat similar to those theorems:

\smallskip
\textbf{Theorem 6.} A rectangle is called good, if one of its sides is integer. Then any rectangle which is
divided into good rectangles is good.
\smallskip

Many nice proofs were invented for this fact. S. Wagon [19] has published a collection of 14 proofs, 
and that collection is far from being complete. In [19] he also explains, that
some of those proofs can be generalized to higher dimensions, and arbitrary additive subgroup of $\RR$
instead of $\ZZ$. Our technique also provides yet another proof for this fact and 
for its generalization:

\smallskip
\textbf{Theorem 7.} Given an additive subgroup $G$ of $\RR$, and number $K \leq N$, an $N$-dimensional
parallelepiped will be called good if lengths of its sides in at least $K$ different directions
belong to $G$. Then any parallelepiped which is divided into good parallelepipeds is good.
\smallskip

Theorems 6 and 7 is useful for some very natural combinatorial riddles.
Once a seven-year-old boy asked his dad why couldn't he fill a $6 \times 6 \times 6$ box 
with $1 \times 2 \times 4$ bricks. His dad happened to be a mathematician and developed
some algebraic theory (N.G. de Bruijn [6]) to answer the question.
Generalization of these question are theorems 6 and 7. More applications of those
theorems to combinatorial riddles will be mentioned after the proofs.

In this paper we prove all the above theorems. Theorems of that kind can have two directions:
to prove that when a certain algebraic condition is satisfied then the parallelepiped 
can be decomposed into parallelepipeds of prescribed kind, and to prove that when the condition
is not satisfied there is no decomposition. The first direction is done by specific 
construction of decomposition (in some cases the existence of decomposition is obvious, 
but not always). The second direction will be performed by construction of
some additive function. The functions will be different for different theorems,
but there are several common points in applying those additive functions.

So, before proving the theorems, we shall explain some general ways of
constructing of additive functions over parallelepipeds.

\section{Additive functions} 

\textbf{Definition.} A function $f$ over parallelepipeds with parallel faces is called additive 
if for any parallelepiped $P$ which is dissected by a plane parallel to its faces
into two parallelepipeds $P_1$ and $P_2$ then \[f(P) = f(P_1)+f(P_2)\]

\textbf{Claim.} For any additive function, and a parallelepiped $P$ subdivided into $n$ 
parallelepipeds $P_1, P_2, ..., P_n$ 
\[ f(P) = f(P_1)+ ... +f(P_n) \] 
We can formulate and prove a more subtle claim.
Assume we have a coordinate system with axes parallel to the edges of given parallelepipeds.
\smallskip

\textbf{Definition.} A function $f$ over some subset $S$ of a set of parallelepipeds 
with parallel faces is called additive 
if for any parallelepiped $P$ which is dissected by a plane parallel to its faces
into two parallelepipeds $P_1$ and $P_2$ such that $P_1, P_2 \in S$ then $f(P) = f(P_1)+f(P_2)$
\smallskip
Let $P$ be a parallelepiped subdivided into $n$ parallelepipeds $P_1, P_2, ..., P_n$ 
Then each of those parallelepipeds can be defined by its bounds in each coordinate.
Assume that $X_k$ is set of all bounds in coordinate $x_k$ for parallelepipeds $P_1, P_2, ..., P_n$.

\smallskip
\textbf{Claim'.} Suppose $f$ is defined and additive on parallelepipeds whose bounds 
in coordinate $x_k$ belong to $X_k$ for all $k$. Then  
\[ f(P) = f(P_1)+ ... +f(P_n) \] 
The motivation to formulate claim' is the following.
As we shall see soon, for some of our theorems on parallelepipeds 
we shall need to construct some $\QQ$-linear functions over $\RR$. 
The construction uses Hamel basis and hence the Axiom of Choice. 
It would be unnatural if some fact about cutting brick into finite number 
of pieces would depend on the Axiom of Choice, and indeed it doesn't.
To avoid using the Axiom of Choice, one can construct the $\QQ$-linear 
functions not on the whole $\RR$, but on its relevant 
finite dimensional over $\QQ$ subspace. 

To keep the ideas transparent, we shall talk about Hamel basis, 
but we want to remark that the same proof works without Zermelo 
as well. Those few readers who refuse to accept 
the Axiom of Choice, will be able to translate all our proofs
into choice-free framework, using the last definition and claim'.

The proof of the claim can be divided into several sub-claims, 
each of those will be simple induction.

\textbf{Claim 1.} Consider a finite family of $k$ parallel planes, 
parallel to a couple of faces of the parallelepiped $P$.
If they subdivide the parallelepiped $P$ into 
parallelepipeds $p_1, p_2, ..., p_n$ then 
\[ f(P) = f(p_1)+ ... +f(p_n) \] 

\textbf{Claim 2.} Consider $q$ finite families of parallel planes, 
each family parallel to a couple of faces of the parallelepiped $P$.
If they subdivide the parallelepiped $P$ into 
parallelepipeds $p_1, p_2, ..., p_n$ then 
\[ f(P) = f(p_1)+ ... +f(p_n) \] 

\textbf{Proof of all claims} 

Claim 1 follows directly from the definition by induction over $k$.
\medskip

Claim 2 follows by induction over $q$. The base of induction, $q=1$
is claim 1. The step of induction is the following.

Take one family of parallel planes, which subdivides
the parallelepiped $P$ into $P_1, ..., P_m$. Each $P_j$ is subdivided 
by only $k-1$ planes into its parts $p_i$ hence by induction
\[ f(P_1)+f(P_2)+...+f(P_m) = f(p_1)+ ... +f(p_n) \] 
while by claim 1 
\[ f(P) = f(P_1)+f(P_2)+...+f(P_m) \] 
Hence
\[ f(P_1)+f(P_2)+...+f(P_m) = f(p_1)+ ... +f(p_n) \] 
QED
\medskip

The claim in the beginning of the section follows from claim 2 
in the following way. Prolong all planes which are faces of the parallelepipeds
of the subdivision. They cut the original parallelepipeds $P_1, ..., P_m$
into smaller parts: $p_1, p_2, ..., p_n$. Therefore, by claim 2 
\[ f(P) = f(p_1)+ ... +f(p_n) = f(P_1)+f(P_2)+...+f(P_m)\] 
QED.

Here we see, that this proof works actually also for the subtler version, claim'.
\medskip

We shall use two constructive ideas to build additive functions over parallelepipeds:

\textbf{First idea.} Let $\phi_1, \phi_2, ... , \phi_N$ be a set of additive functions
over real variables. Let $a_1, a_2, ... , a_N$ be lengths of sides in directions $x_1, x_2, ... , x_N$ 
of parallelepiped $P$. Define
\[ f(P) = \phi_1(a_1) \cdot \phi_2(a_2) \cdot ... \cdot \phi_N(a_N) \]
Then it is clear that this function is additive. Any linear combination of those 
and any polylinear function in $a_1, a_2, ... , a_N$ is additive as well.

\textbf{Second idea.} Parallelepiped $P$ has $2N$ faces: a lower bound in each coordinate
and a higher bound in each coordinate. The faces that correspond to lower bounds 
will be called lower faces, the other faces will be called upper faces.
A vertex of $P$ will be called a black vertex of $P$ if it is contained in
even number of lower faces, otherwise it will be called a white vertex of $P$.

Suppose we have any function on $F:\RR^N \rightarrow \RR$.
Then we can define an additive function of parallelepipeds: 

$f(P) = $ sum of $F$ over black vertices of $P$ minus sum of $F$ over white vertices of $P$.

It is easy to see that the function $f$ is additive over parallelepipeds.

\section{Lower dimensions} 

\textbf{Proof of theorem 1.} Let $w,h$ be width and height of a rectangle.

One direction is obvious: if $\frac{w}{h} = \frac{m}{n}$ 
where $m,n \in \NN$ then $\frac{w}{m} = \frac{h}{n} = a$ 
and then the rectangle can be tiled by the squares with side $a$. 

Let $f(x)$ be a $\QQ$-linear function. Then function $F(rectangle) = w \cdot f(h) - h \cdot f(w)$
is additive, and zero on squares. If width and height of the rectangle is non-commensurable 
we can construct such a $\QQ$-linear function that $f(w) = 0$ , $f(h) = 1$ 
(by choosing a basis of $\RR$ over $\QQ$ which contains both $w$ and $h$).
Then $F(rectangle) \neq 0$.
But if the rectangle would be decomposable into squares, $F(rectangle)$ would be 0
since the function is additive.

Q.E.D.

Another additive function which leads to proof is the following. 

Construct a $\QQ$-linear function $f(x)$ which is positive on $w$ and negative on $h$.
Consider function $F(rectangle) = f(width) \cdot f(height) $. It is nonnegative on each
square and negative on the original square.

\textbf{Another proof of theorem 1.}  (Private communication, A. Kanel-Belov).
Consider the certain sub-division of the certain rectangle into squares. 
Denote $x_i$ the sides of the squares. The subdivision defines certain linear equations over 
the numbers $x_i$, such as: if a certain interval in the picture is presented as sum of sides 
of different subsets of squares, then those sums should be equal. 
Two more equations claim that sum of all sides of squares touching 
the lower rectangle's border is $w$, and sum of two rectangle's squares 
touching the right rectangle's border is $h$.
All coefficients in that system of equations are rational, 
except for $w$ and $h$. After applying scaling we shall assume $h = 1$, 
and then $w$ will be the aspect ratio.
The configuration of squares solves the problem for given $w$ if and only 
if the system of linear equations has a solution in nonnegative real numbers. 
When we apply Gauss method to solve the system of all equations 
except the one containing $w$, we shall get either a single rational 
solutions (because coefficients are rational), 
or an infinite family of solutions, which depends linearly 
(with rational coefficients) upon a finite number of parameters. 
A solution we get from Gauss methods is limited by several inequalities, 
corresponding to non-negativity of all $x_i$. 
So, if we have an infinite family of solution 
(and that is the only way to get irrational $w$) 
then $w$ moves in certain limits, between two rational limiting values, 
$a$ and $b$.
In such a case, we might write $x_i = k_it + m_i$, and for every $t$ 
we shall get the same configuration of squares, but of different sizes. 
Say at $t = 0$ we shall get rectangle of width $a$, at $t = 1$ we shall 
get a rectangle of width $b$, and at some intermediate value we shall get width $w$. 
Sides of all squares change as a linear function of $t$. 
So their areas are convex quadratic functions, 
unless they are unchanged. But height is constantly 1, and w changes linearly, 
hence the area changes linearly. That's a contradiction. 
Hence, there are no configurations which give infinite families 
of solutions, only those,  that give a single rational solution. Q. E. D.

\textbf{Proof of theorem 2.} (About 3-dimensional parallelepipeds and bars.)

Unlike theorem 1, here both directions are non-obvious.

First direction. Assume that $a,b,c$ are linearly dependent over $\ZZ$.

Then non-trivial linear combination can have no more than 1 zero coefficients.
If it has 1 zero coefficient it has also a positive coefficient and a negative 
coefficient. So, the condition takes form $ka=mb$, where $k,m\in\NN$.
Then one can take two families of planes - one cutting $a$-sides into $m$
equal parts and another cutting $b$ sides into $k$ equal parts, then the
parallelepipeds will be decomposed into bars.

If in the linear combination all coefficients are nonzero, there should be two coefficients 
of one sign and one of different sign. 
So, without loss of generality we can write $ka+mb=nc$ where $k,m,n \geq 0$ and integer.
Therefore the sides $c$ can be divided by a plane into two 
parts : $ka/n$ and $mb/n$.
Both parts have a pair of commensurable sides in different directions, 
and for those cases we have solved the problem already.

Second direction: suppose the parallelepiped is dissected into bars. 
Let $f, g$ be $\QQ$-linear functions over $\RR$. 
Consider a parallelepiped $P$ which is defined by his 3 side lengths : length $a$ in $x$ direction, 
length $b$ in $y$ direction, and length $c$ in $z$ direction.
\[ F(P) = \left| \left( \begin{array}{ccc}
a & b & c \\
f(a) & f(b) & f(c) \\
g(a) & g(b) & g(c) \end{array} \right) \right| \] 
It follows from definition that $F(bar)=0$ for any bar, and $F$ is additive.

If $a, b, c$ are linearly independent over $\QQ$, we can construct a $\QQ$-linear function $f$
such that $f(a)=f(c)=0 ; f(b) = 1$ and a a $\QQ$-linear function $g$
such that $g(a)=g(b)=0 ; g(c) = 1$. Then $F(P)$ is non-zero and hence it cannot be dissected 
into bars, contradiction, Q. E. D.

\textbf{Proof of theorem 3.} Here we try to decompose 
parallelepipeds of type $a \times a \times b \times b$
into parallelepipeds of type $a \times a \times a \times b$.
If $a,b$ are commensurable it can be decomposed even into cubes. 

Let $P$ be a 4-dimensional parallelepiped $a \times b \times c \times d$
\[ F(P)= 
\left| \begin{array}{cc}
a & c\\
f(a) & f(c)\end{array}\right| \cdot 
\left| \begin{array}{cc}
b & d\\
f(b) & f(d)\end{array}\right| \] 
Where $f$ is a $\QQ$-linear function over $\RR$. $F$ is polylinear in $a,b,c,d$ hence additive.

$F$ is zero on parallelepipeds of types 
$x \times x \times x \times y$, $x \times x \times y \times x$ , 
$x \times y \times x \times x$, $y \times x \times x \times x$. 

If $a, b$ are non-commensurable we can choose $f$ such that
that $f(a)=0,f(b)=1$. Then
\[F(a,a,b,b)=\left| \begin{array}{cc}
a & b\\
f(a) & f(b)\end{array}\right|^2 = a^2 \neq 0\]
Hence it is not decomposable into parts, on which $F$ is zero.

\textbf{Proof of theorem 4.} Here we try to decompose 
parallelepipeds of type $a \times a \times a \times b$ 
into parallelepipeds of type $x \times x \times y \times y$ .

Suppose, such division is possible and let 
\[ F(cuboid)=f(x)f(y)f(z)f(t)\] 

Then, $F$ is $\QQ$-linear in each variable and $F(x,x,y,y)=F(x,y,y,x)=F(x,y,x,y)=f^2(x)f^2(y)\geq0$. 
Now let $f$ be such $\QQ$-linear function, that $f(a)=-1,f(b)=1$. We have that $F$ is additive. 
But $F$ is nonnegative on all cuboids of type $x \times x \times y \times y$ and $F(a,a,a,b)=-1<0$.

Contradiction.

\section{Positive basis and higher dimensions} 

\textbf{Positive basis theorem.} Consider vectors $v_1,...,v_n \in \QQ^k$, in a half space defined by linear functional $l:\QQ^k \longrightarrow \RR$:
for all $j$, $l(v_j) > 0$. 
Then there exists a basis of $k$ vectors, $e_1, ...,e_k$, lying in the same half space defined by $l$ and having all rational coordinates, 
such that vectors $v_1,...,v_n$ are linear combinations of $e_1, ...,e_k$ with nonnegative coefficients.

\textbf{Conclusion.} If positive real numbers $x_1,...,x_n$ span $k$-dimensional linear space over $\QQ$, one can find positive numbers $e_1, ...,e_k$ such that $x_1,...,x_n$ are linear combinations
of $e_1, ...,e_k$ with nonnegative rational coefficients (or even nonnegative integer coefficients).

The conclusion is just a special case of the positive basis theorem - 
choose largest possible linearly independent subset 
of $\{x_1,...,x_n\}$, it spans a $k$-dimensional linear space over $\QQ$, 
containing each $x_j$, and positive number form a half-space of that
subset, hence there are elements with that property. When You found a 
positive basis over $\QQ$, just divide each element of the basis by common
denominator of all its coefficients.
\smallskip
The positive basis theorem (or rather the conclusion) will be used to prove 
one direction of theorem 5, which is perhaps the most interesting theorem
in the sequence, so we shall prove it before we get to prove theorem 5.
\smallskip

\textbf{Proof of positive basis theorem.} Without loss of generality, we may assume 
that $l(v) = \left\langle C,v \right\rangle $ where $\left\langle C,C \right\rangle = 1$

$\QQ^k$ is dense in  $\RR^k$.
Therefore rays generated by vectors in $\QQ^k$ cut the hyperplane $H=\{v|l(v)=1\}$
over a dense subset.

Consider a $k-1$-dimensional regular simplex in $H$, with vertices $w_1, w_2, ..., w_k$ and 
center $C$, such that the distance from $C$ to edges of simplex is $d$.
Then any vector $V$ in $H$ such that the angle between $ C,V $  is less than $arctan(d)$ 
is in convex hull of $w_1, w_2, ..., w_k$. Hence any vector such that 
the angle between $ C,V $  is less than $arctan(d)$ is linear combination 
of $w_1, w_2, ..., w_k$ with positive coefficients.
The distance between hyperplane in $H$ formed by points $x_1, x_2, ...,x_{k-1}$ and $C$ 
depend continuously on $x_1, x_2, ...,x_{k-1}$, if points are far enough from each other.
So, for any small $\epsilon$ there is $\delta$ such that if $|w_j-u_j|< \delta$, where $u_j \in H$
then the distance from $C$ to all faces of the simplex $u_1, u_2, ..., u_k > d-\epsilon$.
Then any vector $V$, such that tangence of the angle between $V,C$ is less than $d-\epsilon$
is a linear combination of $u_1, u_2, ..., u_k$ with positive coefficients.

Now choose $d$ and $\epsilon$ so that $d-\epsilon$ will be strictly bigger then tangence
of angle between $C$ and $v_i$ for all $i=1, 2, ..., n$.
Choose $u_j$ sufficiently close to $w_j$, so that $u_j$ would be positive multiple of 
some rational vector $e_j$. Obviously, $e_j$ will be in the correct half-space and 
vectors $v_1, v_2, ..., v_n$ will be positive linear combinations of $u_1, u_2, ..., u_k$
and hence of $e_1, e_2, ..., e_k$. 

Q. E. D.
\smallskip

\textbf{Proof of theorem 5.} We try to subdivide the $n$-dimensional parallelepiped, 
whose sides in $n$ different directions are $a_1, a_2, ..., a_n$ into $k$-bars.

Suppose the span of $a_1, a_2, ..., a_n$ over $\QQ$ of the sides of the parallelepiped 
has dimension $k$ or less.
Then, by the conclusion from the positive basis theorem, there are such positive 
$e_1, e_2, ..., e_k$ that $a_1, a_2, ..., a_n$ are their linear combinations
with nonnegative integer coefficients. 

So, we can build $n$ families of parallel planes, each family will split side
in direction $j$ of length $a_j$ into intervals of lengths $e_1, e_2, ..., e_k$.
Therefore, the parallelepiped will be subdivided into parallelepipeds, whose 
all sides are $e_1, e_2, ..., e_k$, so they have no more than $k$ different sides.

Suppose the span of $a_1, a_2, ..., a_n$ over $\QQ$ of the sides of the parallelepiped 
has dimension greater than $k$ . Then we can choose a subset of $k+1$ linearly independent
over $\QQ$ numbers out of $a_1, a_2, ..., a_n$. We may assume without loss of generality 
that those are $a_1, a_2, ..., a_{k+1}$.
Construct $\QQ$-linear functions  $f_1, f_2, ..., f_{k+1}$ from $\RR$ to itself,
such that for $i,j \leq k+1$ 
\[f_i(a_j) = \delta_{i,j}\]
For any parallelepiped $P$ whose sides in $n$ directions are $l_1, l_2, ..., l_n$
in this order define a  $(k+1) \times (k+1)$ matrix $M(P)$ with entries $m_{i,j} = f_i(l_j)$
\[F(P) = det(M(P)) \cdot l_{k+2}l_{k+3}\cdot...\cdot l_n \]
It is obvious that $F$ is additive (since it is polylinear over $l_j$), 
that it is 0 on bars (since two columns of the matrix are equal) and that it is 
non-zero over the original parallelepiped.

Therefore the original parallelepiped is not decomposable into bars. Q.E.D.

\section{"Good" rectangles and parallelepipeds} 

Now we shall prove a claim which is slightly more general than theorem 6,
using the second idea for constructing additive functions.

\textbf{Definition} Let $G$ be a fixed additive subgroup of $\RR$ ($\ZZ$ is only one possible example).
An $n$-dimensional  parallelepiped will be called good, if at least 
one of its sides belongs to $G$. 

\textbf{Theorem 6'}. A parallelepiped, which is partitioned into good parallelepipeds,
is good.

\smallskip
\textbf{Proof of theorems 6 and 6'} We choose coordinates in $\RR^n$ such that the axes 
go along the sides of the parallelepiped. $G^n$ and all its shifts by vectors
(which are elements of the factorgroup $\RR^n / G^n$) will be called lattices.

Let $F$ a function $\RR^n \rightarrow \RR$ whose value depends only on the lattice
(invariant with respect to shifts by $G^n$). This function is not defined uniquely,
it depends on the choice of a function $\RR^n / G^n \rightarrow \RR$, we shall choose 
it later.

Recall, that a vertex of parallelepiped $P$ is called a black vertex of $P$ if it is contained in
even number of lower faces, otherwise it will be called a white vertex of $P$.

We define an additive function of parallelepipeds: 

$f(P) = $ sum of $F$ over black vertices of $P$ minus sum of $F$ over white vertices of $P$.

It is obvious, that if we have a side length belonging to $G$, all pairs of vertices 
connected by a parallel side cancel out. Therefore, if the original parallelepiped is
splittable into good, then any $f$ of that kind is 0 on it.

But if the original parallelepiped has no sides in $G$, then all its vertices belong 
to different lattices, hence we can require that $F$ would be 1 on one of its black 
vertices and 0 on all its other vertices, black and white. Then $f$ will be 1
on the original parallelepiped. Contradiction.

QED.
\smallskip

For one of the proofs of the last theorem we shall need the following lemma

\smallskip
\textbf{Lemma} 
For any function $\phi:\RR \rightarrow \RR$ denote $\Delta _ {a} \phi (x) = \phi(x+a) - \phi(x) $

Let $x \geq 0$ and $a_1, a_2, ..., a_n > 0$.
Then $\Delta _ {a_1} \Delta _ {a_2} ... \Delta _ {a_n} x^k $ is positive for $n\leq k$ and 0 for $n>k$.
\smallskip

\textbf{Proof}
Assume $\phi(x)$ is smooth, like $x^k$ for example.
It is easy to see that for each $a>0$
\[\Delta _ {a} \phi(x) = \phi(x+a) - \phi(x) = \int^{x+a}_{x} \phi'(x) dt \]
Iterating this formula, we get 
\[\Delta _ {a_1} \Delta _ {a_2} ... \Delta _ {a_n} \phi(x) = \int^{x+a_1}_{x}\int^{x+a_2}_{x} ... \int^{x+a_n}_{x} \phi^{(n)}(x) dt \]
Here $\phi^{(n)}$ is $n$-th derivative. By applying Lagrange theorem, we conclude that there is a point $y$ between $x$ and $x+a_1+...+a_n$, 
such that 
\[\Delta _ {a_1} \Delta _ {a_2} ... \Delta _ {a_n} \phi(x) = \phi^{(n)}(y) \cdot a_1\cdot a_2\cdot ... \cdot a_n\]

Substitute $\phi(x) = x^k$ and lemma becomes obvious.

\smallskip
\textbf{Proofs of theorem 7. } 

\smallskip
\textbf{First proof. } Like in proof of theorem 6, consider lattices, which are $G^n$ and all its shifts by vectors.

Let $F$ a function $\RR^n \rightarrow \RR$:
\[ F(x_1, x_2, ...,x_n) = \alpha(x_1, x_2, ..., x_n)\cdot (x_1 + x_2 + ... + x_n)^{k-1}  \]
Here $\alpha(x_1, x_2, ..., x_n)$ is a function, whose value depends only on the lattice
(invariant with respect to shifts by $G^n$). 

Now generate function over parallelepipeds $f(P) = $ sum of $F$ over black vertices 
of $P$ minus sum of $F$ over white vertices of $P$.
	
By the above lemma, it is 0 on all good parallelepipeds.

If the original parallelepiped is not good, then each lattice contains 
no more than $2^{k-1}$ of its vertices. Therefore, sum with signs of its vertices 
belonging to a certain lattice, containing one of the vertices, is $\alpha$ 
of that lattice times a non-zero number (again, by the lemma).
Therefore, we can choose $\alpha$ so that $f$ over original parallelepiped
would not be 0, hence it is not decomposable into good parallelepipeds.

\smallskip
\textbf{Second proof. } Suppose less than $k$ sides in different directions of 
the parallelepiped belong to $G$, however, a partition into good parts exists. 
Denote the sides in different directions $a_1, a_2, ..., a_n$, so that only 
the first of those belong to $G$, so the last $a_{k+1}, ..., a_n$ don't belong to $G$.

Consider $n-k+1$-dimensional face of the parallelepiped with sides $a_{k+1}, ..., a_n$.
Our partition generates a partition of this face into good parallelepipeds 
in the sense of theorem 6', and it leads to a contradiction, Q. E. D.
\smallskip

\textbf{A few examples: }

\textbf{Example 1. } A question of seven-year-old F. W. de Bruijn to his dad [6]: is it possible to fill  
a $6 \times 6 \times 6$ box with $1 \times 2 \times 4$ bricks?

The answer is no. A brick will be called good if one of its sides is integer multiple of 4.
Each brick is good, the box is bad, so by theorem 6' it can't be filled.

\textbf{Example 2. } 24th Tournament of Towns, spring, junior group, training version, problem 5.
Is it possible to tile $2003 \times 2003$ board by vertical $1 \times 3$ rectangles and horizontal
$2 \times 1$ rectangles?

The answer is no. Contract the board and the tiles by factor 3 in vertical direction and by factor
2 in horizontal direction. All tile will have one integer side, so by theorem 6 the board should 
have an integer side, if we can tile it. But it doesn't.

\textbf{Example 3. } 26th Tournament of Towns, autumn, senior group, main version, problem 5.
Let $A$ and $B$ be rectangles. Show that if one can compose a rectangle similar to $B$ 
from rectangles equal to $A$, then also vice versa, one can compose a rectangle similar to $A$ 
from rectangles equal to $B$. 

Proof. Assume that $A$ has non-commensurable sides $w$ and $h$.
A rectangle $B'$ similar to $B$ can be composed from rectangles equal to $A$.
If we define a good rectangle as a rectangle having a side which is integer multiple of $w$,
we see by theorem 6 that $B'$ has a side which is $nw$, where $n$ is integer.
In the same way, we see that it has a side which is $mh$, where $m$ is also integer.
Those two sides are different, since $w, h$ are non-commensurable. 
So, gluing $m \times n$ copies of equally-oriented versions of $B$, we shall get a rectangle 
similar to $A$.
If sides of $A$ are commensurable, than sides of any rectangle composed of its equal 
copies are obviously commensurable, and the statement is obvious. 

Of course, for all those examples there are more elementary proofs, but those require
certain inventiveness, while theorems 6 and 6' create a general and obvious approach.

\medskip
To conclude: constructing a crafty additive function 
can give an elegant solution even for a hard problem.

\end{document}